\documentstyle[12pt,amstex,amssymb]{article}
\newcommand{\Aa}{\Bbb{A}}
\newcommand{\C}{\Bbb{C}}

\newcommand{\R}{\Bbb{R}}
\newcommand{\Z}{\Bbb{Z}}

\newcommand{\Q}{\Bbb{Q}}
\newcommand{\E}{\Bbb{E}}

\newcommand{\V}{{\Bbb V}}
\newcommand{\W}{{\Bbb W}}

\newcommand{\Ch}{{\cal C}}

\newcommand{\Eh}{{\cal E}}
\newcommand{\Fh}{{\cal F}}

\newcommand{\Hh}{{\cal H}}
\newcommand{\Lh}{{\cal L}}
\newcommand{\Mh}{{\cal M}}
\newcommand{\Oh}{{\cal O}}

\newcommand{\Th}{{\cal T}}

\newcommand{\eB}{\frak{B}}

\newcommand{\eX}{\frak{X}}

\newcommand{\esl}{\frak{sl}}

\newcommand{\ev}{\frak{v}}
\newcommand{\ew}{\frak{w}}

\newcommand{\id}{\operatorname{ id}}
\newcommand{\ab}{\operatorname{ab}}

\newcommand{\ad}{\operatorname{ ad}}

\newcommand{\Div}{\operatorname{Div}}

\newcommand{\RRe}{\operatorname{ Re}\,}
\newcommand{\Imm}{\operatorname{ Im}\,}
\newcommand{\imm}{\operatorname{ im}\,}

\newcommand{\For}{\operatorname{ For}}

\newcommand{\coker}{\operatorname{coker}}
\newcommand{\conn}{\operatorname{conn}}
\newcommand{\pr}{pr}
\newcommand{\pure}{\operatorname{pure}}
\newcommand{\Lie}{\operatorname{ Lie}}
\newcommand{\Li}{\operatorname{ Li}}

\newcommand{\Gr}{\operatorname{ Gr}}

\newcommand{\Ext}{\operatorname{ Ext}}

\newcommand{\pol}{pol}
\newcommand{\GL}{\operatorname{GL}}
\newcommand{\SL}{SL}
\newcommand{\Si}{Si}

\newcommand{\Var}{\operatorname{ Var}}

\newcommand{\Spec}{\operatorname{ Spec}}

\newcommand{\tors}{\operatorname{ tors}}

\newcommand{\diff}{\operatorname{ diff}}

\renewcommand{\mod}{\operatorname{ mod}\,}

\newcommand{\dis}{\displaystyle}

\newcommand{\verk}{\raisebox{0.07cm}{\mbox{\scriptsize $\,\circ\,$}}}

\newcommand{\btd}{\triangledown}
\newcommand{\halb}{\frac{1}{2}}
\newcommand{\silo}{\stackrel{\sim}{\longrightarrow}}

\newcommand{\tei}{\, | \,}

\newcommand{\sto}{\mbox{$-\!\to$}}

\newcommand{\limr}[1]{\dis \lim_{\stackrellow{\longrightarrow}{#1}}}

\def\stackrellow#1#2{\mathrel{\mathop{#1}\limits_{#2}}}
\newcommand{\beweisende}{\hspace*{\fill} q.e.d.}
\newtheorem{theorem}{Theorem}[section]
\newtheorem{punkt}[theorem]{$\!\!\!\,$}
\begin{document}
\thispagestyle{empty}
\vspace*{5cm}
\begin{center}
{\Large\sc  Variations of Hodge--de Rham Structure\\
\bigskip
and Elliptic Modular Units}\footnotemark

\footnotetext{Math. Subj. Classification 
numbers: 11 G 16, 14 D 07 (11 F 03, 14 H 52, \\
19 F 27, 32 G 20).}

\bigskip
\bigskip
\bigskip
\bigskip
\bigskip

{\sc J\"org Wildeshaus \\
\bigskip
\bigskip

Math. Institut der WWU M\"unster\\
Einsteinstr. 62\\
D--48 149 M\"unster
 
\bigskip

wildesh at math.uni-muenster.de
}

\end{center}
\newpage
\setcounter{page}{1}

\section*{Introduction}

According to the general motivic folklore, one expects 
the group of one-extensions of $\Q(0)$ by $\Q(1)$ 
in the category $\Mh\Mh_{\Q} (B)$ of 
{\it mixed motivic sheaves} on a scheme $B$ to be given by
\[
\Ext^1_{\Mh\Mh_{\Q} (B)} (\Q (0) , \Q (1)) = \Gamma (B , 
\Oh^{\ast}_B) \otimes_{\Z} \Q \; .
\]
If $B$ is a smooth and separated scheme over a field 
embeddable into
$\C$, then it is possible to define, as a first approximation to 
$\Mh\Mh_{\Q} (B)$, a smooth sheafified variant $MS^s_{\Q} (B)$ of {\it mixed realizations} \`a la Deligne--Jannsen, and a functorial monomorphism
\[
\Gamma (B , \Oh^{\ast}_B) \otimes_{\Z} \Q \longrightarrow 
\Ext^1_{MS^s_{\Q} (B)} (\Q (0) , \Q (1)) \; .
\]
Its cokernel is enormous. Its image is expected to consist of the 
{\it geometrically motivated} one-extensions. 

As far as $\Ext^1 (\Q (0) , \Q (1))$ is concerned, 
it turns out that a much more 
precise approximation to $\Mh\Mh_{\Q} (B)$ is 
obtained by actually {\it forgetting} part of the data of objects of
$\Mh\Mh_{\Q} (B)$: the category 
$HDR^s_{\Q} (B)$ of {\it variations of $\Q$-Hodge--de Rham structure} 
comes about by leaving away the $l$-adic 
components of $MS^s_{\Q} (B)$. 

A variant of the definition used here, in the case
when $B$ is a point, appeared already in \cite{Ha}, section 1. There,
the interested reader also finds a detailed account of Hodge--de Rham
structures in the context of motives, and their $L$-functions.

The first aim of the present article is to 
popularize the variational point of view, and to illustrate the flexibility
of the resulting formalism by a concrete example. 
%

For an elliptic curve $\Eh$ over $B$ with zero section $i$, we 
denote, letting $\widetilde{\Eh} := \Eh - i (B)$, by $\Lh (\Eh)$ the 
$\Q$-vector space with basis
$( \{ s \} \tei s \in \widetilde{\Eh} (B))$.
Furthermore, define
\[
d : \Lh (\Eh) \longrightarrow \Eh (B) \otimes_{\Z} \Eh (B) 
\otimes_{\Z} \Q \; , \quad
\{ s \} \longmapsto s \otimes s \; .
\]
In \cite{W3}, we constructed, using the so-called {\it 
polylogarithmic extension} on $\widetilde{\Eh}$, a homomorphism
\[
\varphi : \ker (d) \longrightarrow \Ext^1_{HDR^s_{\Q}} (\Q (0) , \Q 
(1))
\]
whose image consists of extensions one will consider as 
geometrically motivated.

The main result of this work shows that $\varphi$ factors 
through
the natural monomorphism
\[
\kappa_B : \Gamma (B , \Oh^{\ast}_B) \otimes_{\Z} \Q 
\longrightarrow \Ext^1_{HDR^s_{\Q} (B)} (\Q (0) , \Q (1)) \; ,
\]
thus obtaining a proof of the {\it elliptic Zagier 
conjecture} on the lowest step $k=2$. The proof is the logical continuation
of the sheaf-theoretical approach developed in \cite{W3}, sections
2 and 3. It is thus independent of the one 
sketched in \cite{W3}, 1.9 (which makes use of the Poincar\'e line bundle,
but not of the material contained in later sections of loc.\ cit.). 
It is also independent of the proof given in \cite{W4},
3.2, 3.9 (which gives a geometrical rather than sheaf-theoretical
construction). For the precise statement of our main result, we
refer to Theorem \ref{32}. 

The proof of \ref{32} is in three steps.

1. Our main technical tool will be Theorem \ref{24}. 
It states that the cokernel of $\kappa_B$ remains unchanged under pullback morphisms of a quite general type.

2. Consider the following special elements of $\ker (d)$:

o) $\{ s \}$ for a torsion section $s \in \widetilde{\Eh} (B)$.\\
i) $\{ s \} - \{ s-t \}$ for $s,s-t \in \widetilde{\Eh} (B)$ and $t \in 
\Eh (B)_{\tors}$.\\
ii) $\{ s + t \} + \{ s -t \} - 2 \{ s \} - 2 \{ t \}$ for $s,t, 
s+t, s-t \in \widetilde{\Eh} (B)$.

We prove that given $S \in \ker (d)$, then the restriction of 
$S$ to some open dense subscheme $B'$ of $B$ will lie in the 
subspace of $\ker (d_{B'})$ generated by such elements.

Because of Theorem \ref{24}, we need to show $\varphi (S) \in \imm (\kappa_B)$ only for the special expressions o), i), ii).

3. In order to do so, we identify explicitly the one-extensions 
of variations on $B \otimes_k \C$ in terms of holomorphic functions. 
For this we use the formulae of 
\cite{W2}, chapter 3. Depending on which kind of relation they come 
from, we call the resulting functions on $B$ {\it elliptic modular 
units of the zeroeth, first and second kind} respectively.

The explicit description over $\C$ of these functions is in fact 
the second aim of this work. The functions occurring in the image of 
$\varphi$ can safely be expected to be of arithmetic interest: for 
the units of the zeroeth kind, this is well-known: they are 
specialization of the Siegel units studied in \cite{KL} and 
elsewhere. In particular, the classical elliptic units occur in 
this framework.

More recently, elliptic modular units of the first kind appeared in 
Kato's construction of Euler systems for modular curves.

The plan of the paper is as follows:

In section 1, we define the Siegel function, which will turn out to 
generate all the functions in the image of $\varphi$.

Section 2 can be read independently of the rest of this article. 
It contains a self-contained introduction of
the category $HDR^s_{\Q}$ of variations of
Hodge--de Rham structure, which will hopefully turn out to be of 
interest in its own right. 

In section 3, polylogarithms enter. We review the construction of 
one-extensions of \cite{W3} in the case of interest to us 
(\ref{31}--\ref{33}). In 
sections \ref{34}--\ref{37}, we describe, following the treatment of \cite{W2}, 
chapter 3, the Hodge--de Rham incarnation of the polylog in the 
case when the elliptic curve $\Eh$ is the universal object over 
some modular curve. We need to slightly modify the explicit 
description given in loc.\ cit.\ in order to be able to transfer 
easily the methods developed in section 4 of \cite{BD} to the 
elliptic case. The main result of this section is Theorem \ref{311}, 
where we identify the extension of variations underlying $\varphi 
(S)$. In \ref{312}, we compare the 
formula to the ``na{\"\i}ve'' one obtained by averaging the 
Siegel function over the divisor $S$.

Sections 4--6 are concerned with elliptic modular units of the 
zeroeth, first and second kind respectively. We need to show that 
for $S$ of the special type o), i), or ii) above, the holomorphic 
function $\varphi (S)^{MHS}$ of \ref{311} descends to the field of 
definition $k$ of $B$. For o), we are able, thanks to the explicit description
of our functions, to connect to the classical 
theory of Siegel units. For i) and ii), we use \ref{24} to restrict to 
the case of torsion sections, which then follows from case o).

In section 7, we conclude the proof of Theorem \ref{32}.\\

I wish to thank G. Banaszak and W. Gajda for the invitation to 
Pozna\'n in November 1996, during which this paper was written up.

I am indebted to the referee for useful criticism, and to Gabi Weckermann for her excellent \TeX ing.\\

{\bf Notation:} We denote by $\Hh^+$ the complex upper half plane.

\section{The Siegel function}

\begin{punkt} \label{11} \em
We start by defining the following elementary functions on $\C 
\times \Hh^+$:
\begin{eqnarray*}
r_1 : \C \times \Hh^+ \longrightarrow \R & , & (z,\tau) \longmapsto 
\RRe (z) - \frac{\RRe (\tau) \cdot \Imm (z)}{\Imm (\tau)} \; , \\
r_2 : \C \times \Hh^+ \longrightarrow \R & , & (z,\tau) 
\longmapsto - \frac{\Imm (z)}{\Imm (\tau)} \; , \\
c_{\Hh} : \C \times \Hh^+ \longrightarrow \Hh^+ \!\!\!\!\!& , & (z,\tau) 
\longmapsto \tau \; , \\
c_{\C} : \C \times \Hh^+ \longrightarrow \C & , & (z,\tau) 
\longmapsto z \; .
\end{eqnarray*}
So we have the equality
\[
c_{\C} = - r_2 c_{\Hh} + r_1 \; .
\]
Furthermore, we let
\[
q_{\Hh} := \exp (2\pi i c_{\Hh}) \; , \quad q_{\C} := \exp (2\pi 
ic_{\C}) \; .
\]
\end{punkt}

\begin{punkt} \label{12} \em {\bf Definition:} (cmp.\ \cite{Ku}, (2.14).)
The {\it Siegel function}
\[
Si : \C \times \Hh^+ - \{ (z,\tau) \tei z \in \Z \oplus \Z \tau \} 
\longrightarrow \C
\]
is given by
\[
Si := - \exp (\pi i \cdot B_2 (-r_2) c_{\Hh}) \exp (-\pi i \cdot 
r_1 (r_2 + 1)) (1 - q_{\C}) \prod^{\infty}_{n=1} (1 - q^n_{\Hh} 
q_{\C}) (1 - q^n_{\Hh} / q_{\C}) \; .
\]
\end{punkt}
Recall the shape of the second Bernoulli polynomial:
\[
B_2 (X) = X^2 - X + \frac{1}{6} \; .
\]

\begin{punkt} \label{13} \em The proof of the following is left to the reader:

{\bf Lemma:} For $\tau \in \Hh^+$ and $z = -r_2\tau + r_1 \in \C$, we 
have
\begin{eqnarray*}
Si (z+1,\tau) & = & \exp (-\pi i \cdot (r_2 + 1)) Si (z,\tau) \; , \\
Si (z+\tau,\tau) & = & \exp (-\pi i \cdot (r_1 + 1)) Si (z,\tau) \; .
\end{eqnarray*}
\end{punkt}

\section{Variations of Hodge--de Rham structure}

\begin{punkt} \label{21} \em As far as this article is concerned, 
the natural Tannakian category in which the relevant one-extensions 
live is that of variations of Hodge--de Rham structure. 

{\bf Definition:} Let $k$ be a field which is embeddable into $\C$, 
$X / k$ smooth, separated and of finite type, $F \subset \R$ a 
field. $HDR^s_F (X)$, the category of {\it variations of mixed
$F$-Hodge--de Rham structure} on $X$ consists of families
\[
(\V_{DR} , \V_{\sigma} , I_{DR,\sigma} , F_{\sigma} \tei \sigma : k 
\hookrightarrow \C) \; ,
\]
where
\begin{enumerate}
\item [a)] $\V_{DR}$ is a vector bundle on $X$, equipped with a 
flat connection $\btd$ which is regular at infinity in the sense of 
\cite{D}, II, remark following D\'efinition 4.5. Further parts of 
the data are an ascending filtration $W_{\cdot}$ by flat 
subbundles, called the weight filtration, and a descending 
filtration $\Fh^{\cdot}$ by subbundles, the so-called Hodge 
filtration.
\item [b)] $\V_{\sigma}$ is a variation of mixed $F$-Hodge 
structure ($F$--$MHS$) on $X_{\sigma} (\C)$ which is admissible in the 
sense of \cite{Ka}.
\item [c)] Denote by $\For_{\Oh}$ the forgetful functor assigning to a 
variation of $F$--$MHS$ the underlying flat bifiltered vector bundle. 
$I_{DR,\sigma}$ is an isomorphism
\[
\For_{\Oh} (\V_{\sigma}) \silo \V_{DR} \otimes_{k,\sigma} \C
\]
of flat bifiltered vector bundles.
\item [d)] For any $\sigma : k \hookrightarrow \C$, complex 
conjugation defines a diffeomorphism
\[
c_{\sigma} : X_{\sigma} (\C) \silo X_{\overline{\sigma}} (\C) \; .
\]
For a variation of $F$--$MHS$ $\W$ on $X_{\overline{\sigma}} (\C)$, we 
define a variation $c^{\ast}_{\sigma} (\W)$ on $X_{\sigma} (\C)$ as 
follows: the local system and the weight filtration are the 
pullbacks via $c_{\sigma}$ of the local system and the weight 
filtration on $\W$, and the Hodge filtration is the pullback of the 
conjugate of the Hodge filtration on $\W$. The functor 
$c^{\ast}_{\sigma}$ preserves admissibility.

$F_{\sigma}$ is an isomorphism of variations
\[
\V_{\sigma} \silo c^{\ast}_{\sigma} (\V_{\overline{\sigma}})
\]
such that $c^{\ast}_{\overline{\sigma}} (F_{\sigma}) = F^{-
1}_{\overline{\sigma}}$.
\end{enumerate}
\end{punkt}
Furthermore, we require the following: for each $\sigma$, let 
$\iota_{\sigma}$ be the antilinear involution of $\For_{\diff} 
(\V_{\sigma})$, the $C^{\infty}$-bundle underlying $\V_{\sigma}$, 
given by complex conjugation of coefficients. Likewise, let 
$\iota_{DR,\sigma}$ be the antilinear isomorphism
\[
\For_{\diff} (\V_{\sigma}) \silo c^{-1}_{\sigma} (\For_{\diff} 
(\V_{\overline{\sigma}}))
\]
given by complex conjugation of coefficients on the right hand side 
of the isomorphism in c). Our requirement is the formula
\[
\For_{\diff} (F_{\sigma}) = \iota_{DR,\sigma} \verk \iota_{\sigma} = c^{-
1}_{\sigma} (\iota_{\overline{\sigma}}) \verk \iota_{DR,\sigma} \; .
\]

In the category of these data, it is straightforward to define Tate 
twists $F (n)$ for $n \in \Z$: on $F (n)$, the involution 
$F_{\sigma}$ acts by multiplication by $(-1)^n$. The last condition 
we impose is the existence of a system of polarizations: there are 
compatible morphisms
\[
\Gr^W_n \V_{DR} \otimes_{\Oh_X} \Gr^W_n \V_{DR} \longrightarrow 
F_{DR} (-n) \; , \quad n \in \Z
\]
of flat vector bundles on $X$, and polarizations
\[
\Gr^W_n \V_{\sigma} \otimes_F \Gr^W_n \V_{\sigma} \longrightarrow F 
(-n) \, , \quad \sigma : k \hookrightarrow \C \; , \quad n \in \Z
\]
of variations such that the $I_{DR,\sigma}$ and $F_{\sigma}$ and 
the corresponding morphisms for $F (-n)$ form commutative diagrams.

\begin{punkt} \em \label{22} For $k$ and 
$X$ as in {\rm \ref{21}}, we define a map
\[
\kappa_X : \Gamma (X , \Oh^{\ast}_X) \otimes_{\Z} F \longrightarrow 
\Ext^1_{HDR^s_F (X)} (F (0) , F (1))
\]
as follows:

The underlying bifiltered vector bundle is the trivial bundle with 
basis $(e_0 , \frac{1}{2\pi i} \cdot e_1)$, and
\begin{eqnarray*}
\Fh^0 & := & \langle e_0 \rangle_{\Oh_X} \; , \\
W_{-1} & := & \left\langle \frac{1}{2\pi i} \cdot e_1 
\right\rangle_{\Oh_X} \; .
\end{eqnarray*}
For $g \in \Gamma (X , \Oh^{\ast}_X)$, the flat regular connection 
is trivial on $\frac{1}{2\pi i} \cdot e_1$, and maps $e_0$ to
\[
\frac{dg}{g} \cdot \frac{1}{2\pi i} \cdot e_1 \; .
\]
For any embedding $\sigma$ of $k$ into $\C$, the rational structure 
is given by
\[
\left( e_0 - \log g_{\sigma} \cdot \frac{1}{2\pi i} \cdot e_1 , e_1 
\right) \; .
\]
\end{punkt}
\begin{punkt} \label{23} \em Recall the situation in 
the setting of variations of Hodge structure: Let $Z / \C$ be 
smooth, and denote by $\Var_F (Z)$ the category of admissible 
variations of $F$--$MHS$ on $Z (\C)$. By the same construction as in 
{\rm \ref{22}}, we get a map
\[
\kappa^{MHS}_Z = \Gamma (Z , \Oh^{\ast}_Z) \otimes_{\Z} F 
\longrightarrow \Ext^1_{\Var_F (Z)} (F (0) , F (1)) \; .
\]

{\bf Theorem:} $\kappa^{MHS}_Z$ is an isomorphism.

{\bf Proof:} e.g. \cite{W1}, Theorem 3.7. \beweisende
\end{punkt}
From the theorem, we already conclude that the map $\kappa_X$ of 
{\rm \ref{22}} is injective.\\

Let us describe the inverse of $\kappa^{MHS}_Z$: assume given an 
extension

($\ast$) \hspace{-0.75cm} \centerline{$\dis 0 \longrightarrow F (1) 
\longrightarrow \E \longrightarrow F (0) \longrightarrow 0$}

of variations on $Z$. We get an isomorphism of vector bundles
\[
\Fh^0 (\E) \silo F (0) \otimes_F \Oh_X \; ,
\]
hence a splitting of ($\ast$) on the level of bifiltered vector 
bundles. Denote by $e_1$ the base vector ``$2\pi i$'' of the 
constant variation $F (1)$, by $e_0$ {\it the} global section of 
$\Fh^0 (\E)$ mapping to $1 \in F (0)$, and by $\widetilde{e_0}$ 
some multivalued rational flat section of $\E$ mapping to $1$. We have
\[
e_0 - \widetilde{e_0} \in F(1) \otimes_F \Oh_X \; .
\]
Then the theorem tells us that $e_0 - \widetilde{e_0}$ is 
necessarily of the form
\[
\left( \frac{1}{2\pi i} \cdot f \log g + f' \right) \cdot e_1
\]
for some $g \in \Gamma (Z , \Oh^{\ast}_Z) \otimes_{\Z} \Q$ and $f , 
f' \in F$. We have
\[
g \otimes f = \left( \kappa^{MHS}_Z \right)^{-1} (\E) \; .
\]
{\bf Proposition:} Let $k$ be a field which is embeddable into $\C$, $X/k$
smooth, separated and of finite type. For
\[
\E \in \Ext^1_{HDR^s_F (X)} (F (0) , F (1)) \; ,
\]
the following are equivalent:

i) $\E$ lies in the image of $\kappa_X$.

ii) The collection
\[
\For (\E) \in \prod_{\sigma : k \hookrightarrow \C} \Gamma 
(X_{\sigma} , \Oh^{\ast}_{X_{\sigma}}) \otimes_{\Z} F = 
\prod_{\sigma : k \hookrightarrow \C} \Ext^1_{\Var_F (X_{\sigma})} 
(F (0) , F (1))
\]
of extensions of variations underlying $\E$ lies in the image of
\[
\Delta : \Gamma (X , \Oh^{\ast}_X) \otimes_{\Z} F \longrightarrow 
\prod_{\sigma : k \hookrightarrow \C} \Gamma (X_{\sigma} , 
\Oh^{\ast}_{X_{\sigma}}) \otimes_{\Z} F \; .
\]
If i) and ii) are fulfilled, then $\kappa^{-1}_X (\E) = \Delta^{-1} (\For (\E))$.

\begin{punkt} \label{24} \em As for $\coker \kappa_X$, we have the following

{\bf Theorem:} Let $f : X \to Y$ be a morphism of smooth, 
separated $k$-schemes of finite type inducing an isomorphism of the 
schemes of geometrically connected components:
\[
f_{\conn} : X_{\conn} \silo Y_{\conn} \; .
\]
Then the diagram
\[
\begin{array}{ccc}
\Gamma (Y , \Oh^{\ast}_Y) \otimes_{\Z} F & 
\stackrel{f^{\ast}}{\longrightarrow} & \Gamma (X, \Oh^{\ast}_X ) 
\otimes_{\Z} F \\
\kappa_Y \Big\downarrow & & \kappa_X \Big\downarrow \\
\Ext^1_{HDR^s_F (Y)} (F (0) , F (1)) & 
\stackrel{f^{\ast}}{\longrightarrow} & \Ext^1_{HDR^s_F (X)} (F (0) 
, F (1))
\end{array}
\]
is cartesian. In other words, an extension $\E \in \Ext^1_{HDR^s_F (Y)} (F (0) , F (1))$
lies in the image of $\kappa_Y$ if and only if $f^{\ast} \E$ lies 
in the image of $\kappa_X$.

{\bf Proof:} By \cite{W2}, Proposition 3.35, we have
\[
\coker \kappa_Y = \coker \kappa_{Y_{\conn}} \; , \quad \coker 
\kappa_X = \coker \kappa_{X_{\conn}} \; .
\]
Hence our claim follows from the assumption $f : X_{\conn} \silo Y_{\conn}$
and the snake lemma. \beweisende
\end{punkt}
%
%
%
\begin{punkt}\label{26} \em We conclude with another interpretation 
of the isomorphism $\kappa^{MHS}_Z$ of {\rm \ref{23}}.

By \cite{CKS}, Theorem 2.13, for any variation of Hodge structure 
$\V$, there is a unique decomposition of the bifiltered 
$C^{\infty}$-bundle $\V^{\infty}$ underlying $\V$,
\[
\V^{\infty} = \bigoplus_{p,q} \Hh^{p,q} \; ,
\]
such that

i) $\dis W_k \V^{\infty} = \bigoplus_{p+q\le k} \Hh^{p,q}$,\\
ii) $\dis \Fh^p \V^{\infty} = \bigoplus_{p' \ge p} \Hh^{p',q}$,\\
iii) $\dis \overline{\Hh^{p,q}} = \Hh^{p,q} \mod \bigoplus_{p' <q 
\atop q'<q} \Hh^{p',q'}$.\\
It is easily checked that this decomposition is compatible with the 
tensor structure of the category of variations.
\end{punkt}
We thus get a functorial isomorphism
\[
\Theta_{\V} : \V^{\infty} = \bigoplus_{p,q} \Hh^{p,q} \silo 
(\Gr^W_{\cdot} \V)^{\infty} \; ,
\]
which is compatible with tensor products and formation of duals, 
and which satisfies
\[
\Gr^W_{\cdot} \Theta_{\V} = \id_{(\Gr^W_{\cdot} \V)^{\infty}} \; .
\]
\begin{punkt}\label{27}\em
Now assume given a smooth scheme $Z$ over $\C$, and an exact sequence
\[
0 \longrightarrow F (1) \longrightarrow \V \longrightarrow F (0) 
\longrightarrow 1
\]
of admissible variations of $F$--$MHS$ on $Z$. We claim that there is 
a close connection between $\Theta_{\V}$ and the class of $\V$ in
\[
\Gamma (Z , \Oh^{\ast}_Z) \otimes_{\Z} F \stackrellow{=}{{\rm 
\ref{23}}} \Ext^1_{\Var_F (Z)} (F (0) , F (1)) \; .
\]
In this case, the decomposition is already uniquely characterized by 
axioms i) and ii), and exists on the level of holomorphic 
bundles. We have
\begin{eqnarray*}
\Hh^{0,0} & = & \Fh^0 \V \; ,\\
\Hh^{-1,-1} & = & W_{-2} \V \; ,
\end{eqnarray*}
and the isomorphism $\Theta_{\V}$ can be described by expressing 
the images of a basis of rational flat multivalued sections of $\V$ 
in the basis of $\Q (0) \oplus \Q (1)$ given by $e_0$ and $e_1$. The 
result is a matrix of the shape
\[
\left( \begin{array}{cc}
1 & 0 \\ \ast & 1
\end{array} \right) \; .
\]
By Theorem {\rm \ref{23}}, the $\ast$ is of the shape
\[
- \frac{1}{2\pi i} \cdot f \log g \; ,
\]
and we have $g \otimes f = (\kappa^{MHS}_Z)^{-1} (\V)$.
\end{punkt}

\section{Elliptic polylogarithms}
\begin{punkt}\label{31} \em In \cite{W3}, {\rm \ref{31}}, we 
axiomatized some formal properties of stacks $\Th$ on certain 
schemes, which allowed to construct one-extensions in $\Th (B)$ 
from specific linear combinations of symbols on the
Mordell--Weil group $\Eh (B)$ of an elliptic curve $\Eh$ over $B$.

Let $\Ch$ denote the category of schemes $B$ which are smooth, 
separated and of finite type over some field, which is embeddable into $\C$. 
Then by loc.\ cit., {\rm \ref{32}} c),
\[
B \longmapsto HDR^s_{\Q} (B)
\]
is such a stack on $\Ch$. In particular, we have, using the 
notation of loc.\ cit., {\rm \ref{31}}:

(C) For any elliptic curve $\pi : \Eh \longrightarrow B$,
there is given an object of rank two, $R^1 \pi_{\ast} \Q$ in 
$HDR^s_{\Q} (B)$. The formation of $R^1 \pi_{\ast} \Q$ is 
compatible with base change. Write
\[
V_2 := R^1 \pi_{\ast} \Q(1) \; ,
\]
and use the same symbol for the pullback to $\Eh$, or to the complement $\widetilde{\Eh}$ of the zero section.

(E) For any elliptic curve $\Eh / B$,
there is given an {\it Abel--Jacobi map}
\[
[\;] : \Eh (B) \otimes_{\Z} \Q \longrightarrow \Ext^1_{HDR^s_{\Q} 
(B)} (\Q (0) , V_2)
\]
which is compatible with base change. 

(F) ($N=2$ in loc.\ cit.) Consider
\[
[\Delta] \in \Ext^1_{HDR^s_{\Q} (\Eh)} (F (0), V_2)
\]
as a variation of
Hodge--de Rham structure on $\Eh$.
There is given an extension
$\pol^2$ in $HDR^s_{\Q} (\widetilde{\Eh})$ 
of $V_2$ by 
$[\Delta] (1) \, |_{\widetilde{\Eh}}$, the {\it (small) 
polylogarithmic extension}, such that
\begin{eqnarray*}
\pol^1:= \pol^2 / V_2(1) 
& \in & \Ext^1_{HDR^s_{\Q} (\widetilde{\Eh})} (V_2 ,  
\Q (1)) \\
& = & \Ext^1_{HDR^s_{\Q} (\widetilde{\Eh})} (\Q (0) , V_2^{\vee} (1)) \\
& = & \Ext^1_{HDR^s_{\Q} (\widetilde{\Eh})} (\Q (0) , V_2)
\end{eqnarray*}
equals the restriction of $[\Delta]$ to $\widetilde{\Eh}$. Here, the last 
equality is induced by the isomorphism
\[
V_2 \silo V_2^{\vee} (1)
\]
coming from Poincar\'e duality.
\end{punkt}
Furthermore, a certain {\it norm compatibility} (loc.\ cit.,
{\rm \ref{31}} (G)) is satisfied.

Let us remark that the extension $\pol^2$ is unique if one requires as in
loc.\ cit.\ that it be part of a whole projective system 
$(\pol^N)_N$ of extensions.

Note that in loc.\ cit., {\rm \ref{32}} c), we restricted our 
attention to the smaller category $\Ch'$ of schemes which are 
smooth and quasi-projective over some number field. This assumption 
was not used in the proof of (E). As for (F) and norm 
compatibility, we note that since everything is supposed to be 
compatible with change of the base $B$, the construction in (F), as 
well as norm compatibility carry over to the general case because 
the relevant moduli spaces of elliptic curves together with 
finitely many sections are smooth and quasi-projective over $\Q$.

\begin{punkt}\label{32} \em For a scheme $B$ which is smooth, separated and of finite type over 
some field of characteristic $0$, and an elliptic curve $\Eh$ over 
$B$, we denote, slightly modifying the notation of \cite{W3}, by 
$\Lh (\Eh)$ the $\Q$-vector space with basis $( \{ s \} \tei s \in 
\widetilde{\Eh} (B))$. Furthermore, define
\[
d = d (\Eh) : \Lh (\Eh) \longrightarrow \Eh (B) \otimes_{\Z} \Eh 
(B) \otimes_{\Z} \Q \; , \quad \{ s \} \longmapsto s \otimes s \; .
\]
The rest of this article will be concerned with the proof of the 
following

{\bf Theorem:} There is a homomorphism
\[
\varphi = \varphi (\Eh) : \ker (d) \longrightarrow \Gamma (B , 
\Oh^{\ast}_B) \otimes_{\Z} \Q
\]
with the following properties:

a) $\varphi$ is functorial with respect to change of the base $B$.

b) $\varphi$ satisfies norm compatibility: for any isogeny
$\psi: \Eh_1 \longrightarrow \Eh_2$,
whose kernel consists of sections of $\Eh_1$, and any
$s_{1,\alpha} \in (\Eh_1 - \ker (\psi)) (B)$, 
$q_{\alpha} \in \Q$:
\[
d \left( \sum_{\alpha} q_{\alpha} \{ \psi (s_{1,\alpha}) \} 
 \right) = 0 
\iff d \left( \sum_{\alpha} q_{\alpha} \sum_{t \in \ker (\psi) 
 (B)} \{ s_{1,\alpha} + t \} \right) = 0 \; .
\]
If this is the case, then the equality
\[
\dis \varphi \left( \sum_{\alpha} q_{\alpha} \{ \psi 
(s_{1,\alpha}) \} \right) 
= \dis \varphi \left( \sum_{\alpha} q_{\alpha} \sum_{t\in\ker 
(\psi) (B)} \{ s_{1,\alpha} +t \} \right)
\]
holds.

c) If $\Eh = E$ is an elliptic curve over the spectrum $B$ of a 
finite field extension $K$ of either $\Q_p$ or $\R$, then for any 
$S = \sum_{\alpha} q_{\alpha} \{ s_{\alpha} \}$ in the kernel of 
$d$, the absolute value of $\varphi (S) \in K^{\ast} \otimes_{\Z} 
\Q$ satisfies
\[
\log \Big\| \varphi \left( \sum_{\alpha} q_{\alpha} \{ s_{\alpha} 
\} \right) \Big\| = \sum_{\alpha} q_{\alpha} \lambda_K (s_{\alpha}) 
\, ,
\]
where $\lambda_K$ equals the local N\'eron height function (see 
e.g. \cite{Si}, VI).

{\bf Remark:} It suffices to prove \ref{32} for schemes $B$ which are smooth,
separated and of finite type over some field, which is embeddable into $\C$.
\end{punkt}

\begin{punkt} \label{33} \em As a first approximation to Theorem 
{\rm \ref{32}}, we recall that due to the axioms of {\rm \ref{31}}, 
one can construct from $\pol^2$ certain one-extensions of $\Q (0)$ 
by $\Q (1)$ in $HDR^s_{\Q} (B)$. We have:

{\bf Theorem:} There is a homomorphism
\[
\varphi = \varphi (\Eh) : \ker (d) \longrightarrow 
\Ext^1_{HDR^s_{\Q} (B)} (\Q (0) , \Q (1))
\]
with the following properties:

a) as in Theorem {\rm \ref{32}} a).

b) as in Theorem {\rm \ref{32}} b).

c) If $\Eh = E$ is an elliptic curve over $\C$, then for any
$S = \sum_{\alpha} q_{\alpha} \{ s_{\alpha} \}$ in $\ker (d)$,
the extension of Hodge structures
\[
\varphi (S)^{MHS} := \varphi (S)_{\sigma = \id : \C \to \C} \in
\Ext^1_{MHS} (\Q (0) , \Q (1))
\stackrellow{=}{{\rm \ref{23}}} \C^{\ast} \otimes_{\Z} \Q
\]
satisfies
\[
\log \Big\| \varphi \left( \sum_{\alpha} q_{\alpha} \{ s_{\alpha} 
\} \right) \Big\| = \sum_{\alpha} q_{\alpha} \lambda_{\C} (s_{\alpha}) 
\; .
\]

{\bf Proof:} a) and b) is \cite{W3}, Corollary 3.5.
c) is \cite{W3}, Theorem 4.2, together with \cite{Si}, VI, Theorem 3.4, and
\cite{L}, chapter 20, \S\,5. \beweisende
\end{punkt}
\begin{punkt}\label{34} \em Let us describe in explicit terms 
the variation of Hodge 
structure underlying $\pol^2$. We follow the treatment of chapter 3 
of \cite{W2}, where the case of the universal elliptic curve over 
some modular curve was considered.

Fix $n \ge 3$, and let
\[
\pi_n : \Eh_n \longrightarrow Y (n)
\]
denote the universal elliptic curve over the modular curve $Y (n)$ 
``of full level $n$''. $Y (n)$ is a smooth affine scheme over $\Q 
\left( e^{\frac{2\pi i}{n}} \right)$. We shall also consider it as 
a scheme over $\Q$. The scheme $Y (n)_{\conn}$ of geometrically 
connected components equals $\Spec \left( \Q \left( e^{\frac{2\pi 
i}{n}} \right) \right)$.

Let us describe $\pi_n (\C) : \Eh_n (\C) \to Y (n) (\C)$, and 
simultaneously connect to the notation of \cite{W2}. We let
\[
L := L_n := \ker (\GL_2 (\hat{\Z}) \longrightarrow \GL_2 (\Z / n 
\Z)) \; ,
\]
$N := 1$;
\[
V_2 (\hat{\Z}) := \left\{ {a \choose b} \tei a,b \in \hat{\Z} 
\right\} \; ,
\]
and $K:= K_{a,1} := V_2 (\hat{\Z}) \rtimes L$, the semidirect 
product with respect to the natural action of $L$ on $V_2 
(\hat{\Z})$. We let
\[
P_2 := P_{2,a} := \left( \begin{array}{ccc}
1 & 0 & 0 \\ \ast & \ast & \ast \\ \ast & \ast & \ast
\end{array} \right) \le \GL_3 \; ,
\]
and consider $K$ as a subgroup of $P_2 (\Aa_f)$. It is open and 
compact. We have a natural action of $P_2 (\R)$ on $\C \times \Hh^+$:
\[
\left( \begin{array}{ccc}
1 & 0 & 0 \\ a & \alpha & \beta \\ b & \gamma & \delta
\end{array} \right) \in P_2 (\R)
\]
acts by sending $(z,\tau)$ to
\[
\left( (\alpha \delta - \beta \gamma) \cdot \frac{z}{\gamma \tau + 
\delta} + \left( -b \frac{\alpha \tau + \beta}{\gamma \tau + 
\delta} + a \right) , \frac{\alpha \tau + \beta}{\gamma \tau + 
\delta} \right) \; .
\]
Then writing $P'_2 := P_2 \cap \SL_3$, we have
\begin{eqnarray*}
Y (n) (\C) & = & \SL_2 (\Z) \setminus (\Hh^+ \times (\GL_2 
(\hat{\Z}) / L )) \; , \\
\Eh_n (\C) & = & P'_2 (\Z) \setminus (\C \times \Hh^+ \times (P_2 
(\hat{\Z}) / K )) \; ,
\end{eqnarray*}
and $\pi_n (\C)$ is induced by the natural projections
\begin{eqnarray*}
c_{\Hh} : \C \times \Hh^+ & \longrightarrow & \Hh^+ \quad 
\mbox{and} \\
P_2 (\hat{\Z}) & \longrightarrow & \GL_2 (\hat{\Z})
\end{eqnarray*}
respectively. As for the connected components of $Y (n) (\C)$ and 
$\Eh_n (\C)$, one defines
\begin{eqnarray*}
\Gamma & := & \SL_2 (\Z) \cap L \le \GL_2 (\Q) \quad \mbox{and} \\
\Lambda & := & \left( \begin{array}{ccc}
1 & 0 & 0 \\
\Z \atop \Z & \multicolumn{2}{c}{\Gamma}
\end{array} \right) \le P_2 (\Q) \; .
\end{eqnarray*}
Note that the determinant induces isomorphisms
\[
P'_2 (\Z) \setminus P_2 (\hat{\Z}) / K \silo \SL_2 (\Z) \setminus 
\GL_2 (\hat{\Z}) / L \silo (\Z / n \Z)^{\ast} \; .
\]
Choose a set of representatives $R \subset \GL_2 
(\hat{\Z})$ for
$\SL_2 (\Z) \setminus \GL_2 (\hat{\Z}) / L$, and write
\[
p_f := \left( \begin{array}{ccc}
1 & 0 & 0 \\
0 \atop 0 & \multicolumn{2}{c}{g_f}
\end{array} \right) \in P_2 (\hat{\Z})
\]
for any $g_f \in R$. Then we have:
\begin{eqnarray*}
\Eh_n (\C) & = & \coprod_{g_f \in R} \Lambda \setminus (\C \times 
\Hh^+) \, ,\\
Y (n) (\C) & = & \coprod_{g_f \in R} \Gamma \setminus \Hh^+ \; , 
\end{eqnarray*}
and the inclusion of the connected component indexed by $g_f \in R$ 
into
\begin{eqnarray*}
\Eh_n (\C) & = & P'_n (\Z) \setminus (\C \times \Hh^+ \times (P_2 
(\hat{\Z}) / K)) \quad \mbox{and} \\
Y (n) (\C) & = & \SL_2 (\Z) \setminus (\Hh^+ \times (\GL_2 
(\hat{\Z}) / L)) 
\end{eqnarray*}
respectively is given by assigning to the classes of $(z,\tau) \in 
\C \times \Hh^+$ and $\tau \in \Hh^+$ the classes of
$(z,\tau,p_f)$ and $(\tau , g_f)$,
respectively.
\end{punkt}

\begin{punkt}\label{35} \em We have
$\Gr^W_{\cdot} \pol^2 = V_2 \oplus \Q (1) \oplus V_2 (1)$.
One way to describe the canonical isomorphism $\Theta_{\pol^2}$ of
{\rm \ref{26}} is to express a basis of rational flat 
multivalued sections of $\pol^2$ in the corresponding basis of 
$\Gr^W_{\cdot} \pol^2$.
\end{punkt}
As in \cite{W2}, chapter 3, we 
use the parameterization of any of the connected components $\Eh_n (\C)^0$ of $\Eh_n (\C)$ given by  the universal covering map
\[
\pr = \pr_1 : \C \times \Hh^+ \longrightarrow \Eh_n (\C)^0 \; .
\]
We write $(e_1, e_2)$ for the basis of the homology sheaf 
$V_2$, whose value at $\tau \in \Hh^+$ is given by the lines 
connecting $0$ and $1$, resp.\ $0$ and $-\tau$. In order to distinguish the basis $(e_1 , e_2)$ of 
$V_2 \subset \Gr^W_{\cdot} \pol^2$ from the basis $(2\pi i \cdot 
e_1 , 2\pi i \cdot e_2)$ of $V_2 (1) \subset \Gr^W_{\cdot} \pol^2$, 
we follow the notation of loc.\ cit. and write $(\varepsilon'_1 , 
\varepsilon'_2)$ for the basis of $V_2$. We end up with a basis
\[
\eB := (\varepsilon'_1 , \varepsilon'_2 , 2\pi i \cdot 1 , 2 \pi i 
\cdot e_1 , 2\pi i \cdot e_2)
\]
of flat rational multivalued sections of $\Gr^W_{\cdot} \pol^2$.

\begin{punkt}\label{36} \em A basis of flat rational multivalued 
section of $\pol^2 \subset (\pol^2)^{\infty}$ is given by the 
columns of the matrix $P^W_1$ of \cite{W2}, Lemma 3.13. We first 
define its entries:

{\bf Definition:} The {\it $(0,1)$-th elliptic higher logarithm} is 
defined as
\[
\Li_{0,1} := \frac{1}{2\pi i} \left( \sum^{\infty}_{j=0} \log (1 - 
q^j_{\Hh} q_{\C}) + \sum^{\infty}_{j=1} \log (1 - q^j_{\Hh} / q_{\C}) 
\right) + \halb B_2 (-r_2) c_{\Hh} \; .
\]
{\bf Theorem:} Let $\eB = (\varepsilon'_1 , \varepsilon'_2 , 2 \pi 
i \cdot 1 , 2 \pi i \cdot e_1 , 2 \pi i \cdot e_2)$ be the basis of 
{\rm \ref{35}} of flat rational multivalued sections of 
$\Gr^W_{\cdot} \pol^2$. There is a basis of flat rational 
multivalued sections of $\pol^2$ whose image under 
$\Theta_{\pol^2}$ is described by the columns of the matrix
\[
P := \left( \begin{array}{ccccc}
1 & 0 & 0 & 0 & 0 \\[0.2cm]
0 & 1 & 0 & 0 & 0 \\[0.2cm]
-r_2 & r_1 & 1 & 0 & 0 \\[0.2cm]
\Li_{0,1} - \halb r_1 & -\halb r^2_1 - \frac{1}{12} & -r_1 & 1 & 0 
\\[0.2cm]
\halb r^2_2 + \frac{1}{12} & \Li_{0,1} - r_1r_2 - \halb r_1 + 
\frac{1}{4} & -r_2 & 0 & 1
\end{array} \right) \; \begin{array}{c}
\varepsilon'_1 \\[0.2cm] \varepsilon'_2 \\[0.2cm] 2 \pi i \cdot 1 
\\[0.2cm] 2 \pi i \cdot e_1 \\[0.2cm] 2 \pi i \cdot e_2 \end{array}
\]

{\bf Proof:} This is \cite{W2}, Lemma 3.13. There, we used a 
matrix called $P^W_1$, which is best suited as far as norm compatibility 
is concerned (see loc.\ cit., Corollary 3.16). The matrix $P$ is 
obtained from $P^W_1$ by adding $\halb$ times the third column to the 
first and second columns. We thus get a basis of flat rational 
sections which makes the equality of extensions
$\pol^1 = [s]$ of {\rm \ref{31}}.(F) more transparent. \beweisende
\end{punkt}
\begin{punkt}\label{37} \em The matrix $P$ plays a role 
analogous to the one of the matrix $L (z)$ used in \cite{BD}.\\
Following \cite{W3}, we set
\begin{eqnarray*}
c_1 & := & \left( \begin{array}{ccccc}
0 & 0 \, \vline & & & \\
0 & 0 \, \vline & & & \\ \cline{1-3}
 & \hspace*{0.2cm} \vline & 0 \, \vline & 
 \multicolumn{2}{c}{\raisebox{0.3cm}[-0.3cm]{\Large 0}} \\ \cline{3-5}
 & & 1 \, \vline & 0 & 0 \\
\multicolumn{2}{c}{\raisebox{0.3cm}[-0.3cm]{\Large 0}} & 0 \, 
\vline & 0 & 0
 \end{array} \right) \; , \\
c_2 & := & \left( \begin{array}{ccccc}
0 & 0 \, \vline & & & \\
0 & 0 \, \vline & & & \\ \cline{1-3}
 & \hspace*{0.2cm} \vline & 0 \, \vline & 
 \multicolumn{2}{c}{\raisebox{0.3cm}[-0.3cm]{\Large 0}} \\ \cline{3-5}
 & & 0 \, \vline & 0 & 0 \\
\multicolumn{2}{c}{\raisebox{0.3cm}[-0.3cm]{\Large 0}} & 1 \, 
\vline & 0 & 0
 \end{array} \right) \; , \\
d_1 & := & \left( \begin{array}{ccccc}
0 & 0 \, \vline & & & \\
0 & 0 \, \vline & & & \\ \cline{1-3}
1 & 0 \, \vline & 0 \, \vline & 
\multicolumn{2}{c}{\raisebox{0.3cm}[-0.3cm]{\Large 0}} \\ \cline{3-5}
 & & \hspace*{0.2cm} \vline & 0 & 0 \\
\multicolumn{2}{c}{\raisebox{0.3cm}[-0.3cm]{\Large 0}} & 
\hspace*{0.2cm} \, \vline & 0 & 0
 \end{array} \right) \; , \\
d_2 & := & \left( \begin{array}{ccccc}
0 & 0 \, \vline & & & \\
0 & 0 \, \vline & & & \\ \cline{1-3}
0 & 1 \, \vline & 0 \, \vline &
\multicolumn{2}{c}{\raisebox{0.3cm}[-0.3cm]{\Large 0}} \\ \cline{3-5}
& & \hspace*{0.2cm} \vline & 0 & 0 \\
\multicolumn{2}{c}{\raisebox{0.3cm}[-0.3cm]{\Large 0}} & 
\hspace*{0.2cm} \, \vline & 0 & 0
 \end{array} \right) \; .
\end{eqnarray*}
Define $\ev$ as the Lie algebra generated by $c_1 , c_2 , d_1 , d_2$. So
\[
\ev = \left( \begin{array}{ccccc}
0 & 0 \, \vline & & & \\
0 & 0 \, \vline & & & \\ \cline{1-3}
\ast & \ast \, \vline & 0 \, \vline & 
\multicolumn{2}{c}{\raisebox{0.3cm}[-0.3cm]{\Large 0}} \\ \cline{3-5}
\ast & \ast & \ast \, \vline & 0 & 0 \\
\ast & \ast & \ast \, \vline & 0 & 0 \end{array} \right) \subset 
\esl_5 \, ,
\]
a basis being given by
\[
(c_1 , c_2 , (\ad c_2)^l (\ad c_1)^m (d_i) \tei 0 \le l , m \le 1 , 
m+l \le 1) \; .
\]
Inside $\ev$, consider the Lie algebra
\[
\ew := \langle c_1 - d_2 , c_2 + d_1 , (\ad c_2)^l (\ad c_1)^m 
(d_i) \tei m+l = 1 \rangle \; .
\]
 
The Lie algebras $\ew \subset \ev$ correspond to unipotent subgroups 
$W \le V$ of 
$\SL_5$, and $P$ is a multivalued function with values in $W (\C)$. 
Writing $P = (p_{ij})_{1 \le i , j \le 5}$, we have in particular
\begin{eqnarray*}
p_{41} & = & \Li_{0,1} - \halb r_1 \; , \\
p_{52} & = & \Li_{0,1} - r_1 r_2 - \halb r_1 + \frac{1}{4} \; .
\end{eqnarray*}
\end{punkt}

\begin{punkt}\label{38} \em Let $B$ be smooth over $\C$,
$\pi : \Eh \longrightarrow B$ an elliptic curve such that for some $n \ge 3$, the whole
$n$-torsion of $\Eh$ consists of sections of $\pi$. Then we get a 
cartesian diagram
\[
\begin{array}{ccc}
\Eh & \stackrel{f}{\longrightarrow} & \Eh_{n,\C} \\
\pi \Big\downarrow & & \Big\downarrow \pi_{n,\C} \\
B & \stackrel{f}{\longrightarrow} &Y (n)_{\C}
\end{array}
\]
and in particular, a lift of $f$, also denoted
$f : \eX \longrightarrow \coprod_{g_f \in R} (\C \times \Hh^+)$, from the universal cover $\eX$ of $\Eh (\C)$.

Since the formation of $\pol^2$ is compatible with base change, the 
matrix $f^{\ast} P$ describes the isomorphism 
$\Theta_{\pol^2}$. By abuse of notation, we again write $P$ and 
$p_{ij} , 1 \le i,j \le 5$ for $f^{\ast} P$ and its entries.

For any locally closed submanifold $C$ of $\widetilde{\Eh} (\C)$, any 
basis of flat rational sections of $\pol^2 \, |_C$ respecting the 
weight filtration and inducing the basis $\eB$ of $\Gr^W_{\cdot} 
\pol^2 \, |_C$ gives a multivalued function $P'$ satisfying
$P' = P \cdot U$
for a $V (\Q)$-valued function $U$ on $C$, which is therefore 
constant on connected components. Writing $P' = (p'_{ij})_{1 \le i , j \le 5}$,
we have the relations
\begin{eqnarray*}
p'_{31} & = & p'_{53} \; , \\
p'_{32} & = & - p'_{43} \; ,
\end{eqnarray*}
if and only if $U$ is actually $W(\Q)$-valued.
\end{punkt}

{\bf Definition:} (cmp.\ \cite{BD}, 4.1.)
A multivalued function $P'$ with values in $W 
(\C)$ is called a {\it generalized determination} of $P$ if $P$ and 
$P'$ induce the same function with values in
$W (\C) / W (\Q)$.

\begin{punkt}\label{39} \em For the entries of the generalized determination $P' = P \cdot 
U$, we have:

$a_1$) $p'_{32} = - p'_{43} = r_1 + u_1$.\\
$a_2$) $-p'_{31} = -p'_{53} = r_2 + u_2$.\\
$b_1$) $p'_{41} = p_{41} + u_2 r_1 + x$.\\
$b_2$) $p'_{52} = p_{52} -u_1 r_2 + y$.
\end{punkt}
Call a $4$-tuple $(R_1 , R_2 , P_{41} , P_{52})$ of
functions a generalized determination of
$(r_1 , r_2 , p_{41} , p_{52})$
if its entries occur as $p'_{32} = -p'_{43}$, $-p'_{31} 
= -p'_{53}$, $p'_{41}$, and $p'_{52}$ of a generalized determination of 
$P$.

\begin{punkt}\label{310} \em We now imitate the construction of
\cite{BD}, 4.2. With the notation of {\rm \ref{39}}, assume given a 
finite subset $\{ s_{\alpha} \tei \alpha \in I \}$ of $\widetilde{\Eh} 
(B)$, and consider the group
\[
\langle s_{\alpha} \rangle_{\alpha \in I} \le \Eh (B)
\]
generated by the $s_{\alpha}$. Let $\Delta \subset B (\C)$ be a 
simply connected locally closed submanifold, e.g., an open ball. 
For each $s \in \langle s_{\alpha} \rangle_{\alpha \in I}$, choose 
(one-valued!) generalized determinations $R_1$ and $R_2$ of $r_1$ 
and $r_2$ on
\[
s (\Delta) \subset \Eh (\C) \; ,
\]
i.e., functions
$R_i : s (\Delta) \longrightarrow \R$
inducing the same functions modulo $\Q$ as $r_1$ and $r_2$ 
respectively. Furthermore, ensure that these choices are made in a 
way compatible with the group structure, i.e.,
\[
R_i (s (z)) + R_i (t (z)) = R_i ((s+t) (z))
\]
for all $s,t \in \langle s_{\alpha} \rangle_{\alpha \in I}$, and $z 
\in \Delta$.
\end{punkt}
Next, choose for any of the $s_{\alpha} , \alpha \in I$ a 
generalized determination
$(R_1 , R_2 , P_{41} , P_{52})$
of $(r_1 , r_2 , p_{41} , p_{52})$ on
\[
s_{\alpha} (\Delta) \subset \widetilde{\Eh} (C) \; ,
\]
which is compatible with the choices already made.

Observe that any other choice of $R_i$ is of the form $R_i + u_i$ 
for a $\Q$-valued function $u_i$ satisfying
\[
u_i (s (z)) + u_i (t (z)) = u_i ((s+t) (z)) \; .
\]
Different choices therefore lead to the replacements
\begin{eqnarray*}
P_{41} & \longmapsto & P_{41} + u_2 R_1 + x \; , \\
P_{52} & \longmapsto & P_{52} - u_1 R_2 + y \; ,
\end{eqnarray*}
and one concludes:

{\bf Lemma:} If a linear combination $\sum_{\alpha} q_{\alpha} \{ 
s_{\alpha} \}$ satisfies
\[
d \left( \sum_{\alpha} q_{\alpha} \{ s_{\alpha} \} \right) = 
\sum_{\alpha} q_{\alpha} \cdot s_{\alpha} \otimes s_{\alpha} = 0 
\in \Eh (B) \otimes_{\Z} \Eh (B) \otimes_{\Z} \Q \; ,
\]
then the sums
\[
\sum_{\alpha} q_{\alpha} P_{41} \verk s_{\alpha} \quad \mbox{and} 
\quad \sum_{\alpha} q_{\alpha} P_{52} \verk s_{\alpha} \; ,
\]
considered as functions
\[
\Delta \longrightarrow \C / \Q \; ,
\]
are independent of the choice of
\[
(R_1 , R_2 , P_{41} , P_{52}) \; .
\]

\begin{punkt}\label{311} \em With the notation of {\rm \ref{310}}, 
assume given
\[
S = \sum_{\alpha} q_{\alpha} \{ s_{\alpha} \} \in \ker (d) \; .
\]
By the lemma just proved, the functions
\[
\halb \sum_{\alpha} q_{\alpha}  (P_{41} + P_{52}) \verk s_{\alpha}
\]
on open balls in $B (\C)$ glue together to a multivalued function 
on the whole of $B (\C)$, which is well-defined modulo $\Q$. So if 
we define
\[
g_S := \exp \left( - 2\pi i \cdot \halb \sum_{\alpha} q_{\alpha} 
(P_{41} + P_{52}) \verk s_{\alpha} \right) \; ,
\]
then $g_S$ is a multivalued function
\[
B (\C) \longrightarrow \C^{\ast} \; ,
\]
which is well-defined as a function to $\C^{\ast} \otimes_{\Z} \Q$. 
On simply connected open subsets of $B (\C)$, it is representable 
by functions which are differentiable.

From Theorem {\rm \ref{33}}, we recall that $S$ defines an element
\[
\varphi (S) \in \Ext^1_{HDR^s_{\Q} (B)} (\Q (0) , \Q (1)) \; ,
\]
and in particular, an element, denoted by $\varphi (S)^{MHS}$, in
\[
\Ext^1_{\Var_{\Q} (B)} (\Q (0) , \Q (1)) \stackrellow{=}{{\rm 
\ref{23}}} \Gamma (B , \Oh^{\ast} (B)) \otimes_{\Z} \Q \; .
\]
{\bf Theorem:} $\varphi (S)^{MHS}$ and $g_S$ agree as functions from 
$B (\C)$ to $\C^{\ast} \otimes_{\Z} \Q$.

{\bf Remark:} In particular, our local construction of $g_S$ gives 
a holomorphic function, which can be continued to the whole 
universal cover of $B (\C)$, and which represents $\varphi (S)^{MHS}$.

{\bf Proof of Theorem {\rm \ref{311}}:} We imitate the proof of
\cite{BD}, Proposition 4.6. Fix an arbitrary point $b \in B (\C)$, 
and denote by $\omega$ the fibre functor on $\Var_{\Q} (B)$ 
associating to a variation the vector space underlying its fibre at 
$b$. There is a second fibre functor $\omega_0$ on $\Var_{\Q} (B)$ given by
\[
\omega_0 : \V \longmapsto \omega (\Gr^W_{\cdot} \V) \; .
\]
$\omega_0$ and $\omega$ coincide on the subcategory $\Var^{\pure}_{\Q} (B)$ 
of variations with split weight filtration. If $G$ denotes the 
Tannakian dual $G$ of $\Var^{\pure}_{\Q} (B)$, then the dual of 
$\Var_{\Q} (B)$ with respect to $\omega_0$ is a semidirect product
\[
W \rtimes G \; ,
\]
with a pro-unipotent group $W$ (compare \cite{W3}, 2.5).

By \cite{DM}, Theorem 2.13, $\omega$ defines an element in $H^1 (\Q , W 
\rtimes G)$ mapping to zero in $H^1 (\Q , G)$. Since $H^1 (\Q , W) 
= 0$, there is an isomorphism of fibre functors
\[
\omega \silo \omega_0
\]
(see loc.\ cit.). On the other hand, the construction of {\rm 
\ref{26}} gives an isomorphism
\[
\omega_{\C} \silo \omega_{0,\C} \; ,
\]
which is the identity on $\Var^{\pure}_{\Q} (B)$.

Comparing these isomorphisms, we get an automorphism of $\omega_{0,\C}$, 
i.e., an element
\[
w \cdot g \in (W \rtimes G) (\C)
\]
such that $g \in G (\Q)$. We may assume $g = 1$. Then $w \in W 
(\C)$ determines generalized determinants $R_1$ and $R_2$ of $r_1$ 
and $r_2$ on
\[
s (b) \subset \Eh (\C)
\]
for any $s \in \Eh (B)$, which behave additively: for an extension 
$\V$ of $\Q (0)$ by $V_2$, the image of $w$ in $\GL (\omega_0 (\V))$ is 
of the shape
\[
\left( \begin{array}{ccc}
1 & 0 & 0 \\
x & 1 & 0 \\
y & 0 & 1
\end{array} \right) 
\]
where $x$ and $y$ are additive in $\E$, and we set
\[
R_1 (s (b)) := -x \; , \quad R_2 (s (b)) := -y
\]
for $\V = [s]$.

Similarly, if $s \in \widetilde{\Eh} (B)$, then the image $w_s = 
(w_{ij,s})_{1 \le i , j \le 5}$ of $w$ in $\GL (\omega_0 (s^{\ast} 
\pol^2))$ is a generalized determination of the matrix $P (s (b))$.

Now by definition of $\varphi (S)^{MHS}$ (\cite{W3}, 3.3), the 
element $\log (w)$ of $\Lie (W)$ acts on the corresponding 
extension by the matrix
\[
\left( \begin{array}{cc}
1 & 0 \\ \ast & 1
\end{array} \right) \; ,
\]
where the $\ast$ equals
\[
\halb \sum_{\alpha} q_{\alpha} (w_{41,s} + w_{52,s}) \; .
\]
By the recipe given in {\rm \ref{27}}, we have
\[
\varphi (S)^{MHS} (b) = \exp (-2 \pi i \cdot \ast) \; .
\]
\beweisende
\end{punkt}

\begin{punkt}\label{312} \em We remark that we have the equality
\[
\exp \left( - 2 \pi i \cdot \left( \halb (p_{41} + p_{52}) \right) 
\right) = Si^{-1} \; ,
\]
at least up to an eight root of unity.

So if $S = \sum_{\alpha} q_{\alpha} \{ s_{\alpha} \} \in \ker (d)$, then the function $\varphi (S)^{MHS}$ 
differs from the multivalued function
\[
\prod_{\alpha} (Si^{-1} \verk s_{\alpha})^{q_{\alpha}}
\]
by a multivalued function of the shape
\[
\prod_{\alpha} (\exp (2\pi i \cdot F_{\alpha}) \verk 
s_{\alpha})^{q_{\alpha}}
\]
where the $F_{\alpha}$ are polynomials in $\Q [r_1 , r_2]$ of total 
degree smaller or equal to one.
\end{punkt}

\section{Elliptic modular units of the zeroeth kind}

\begin{punkt}\label{41} \em The most visible elements of the kernel of
\[
d : \Lh (\Eh)  \longrightarrow  \Eh (B) \otimes_{\Z} \Eh (B) 
\otimes_{\Z} \Q \; , \quad
\{ s \}  \longmapsto  s \otimes s
\]
are certainly those of the shape
$\{ s \}$,
for a torsion section $s$ of $\Eh$ disjoint from the 
zero section.
\end{punkt}
As in {\rm \ref{38}}, let $B$ be smooth over $\C$, and assume that 
for some $n \ge 3$, the whole $n$-torsion of $\Eh$ consists of 
sections of $\pi$. So we get a cartesian diagram
\[
\begin{array}{ccc}
\Eh & \stackrel{f}{\longrightarrow} & \Eh_{n,\C} \\
\pi \Big\downarrow & & \Big\downarrow \pi_{n,\C} \\
B & \stackrel{f}{\longrightarrow} & Y (n)_{\C}
\end{array} \; .
\]

{\bf Theorem:} Let $s \in \Eh (B)$ be a torsion section. Then
\[
\varphi (\{ s \})^{MHS} = 1 / \Si \verk f \verk s \in \Gamma (B , 
\Oh^{\ast} (B)) \otimes_{\Z} \Q \; .
\]
{\bf Proof:} By {\rm \ref{312}}, both sides differ multiplicatively 
by a function of the shape
$\exp (2 \pi i \cdot F) \verk s$,
where $F$ is a polynomial in $\Q [r_1 , r_2]$. But since $s$ is a 
torsion section, $r_1$ and $r_2$ are rational constants. \beweisende

\begin{punkt}\label{42} \em The functions classically known as {\it 
Siegel units} come about as specializations of $\Si$ to ``torsion 
sections'' of
\[
c_{\Hh} : \C \times \Hh^+ \longrightarrow \Hh^+ \; .
\]

{\bf Definition:} Let $v \in \Q^2 - \Z^2$. The function
$\Si^v : \Hh^+ \longrightarrow \C$
is given by
\[
\Si^v (\tau) := \Si (-v_2 \tau + v_1, \tau) \; .
\]
So $\Si^v$ coincides with the classical Siegel function $g_{(-v_2 , 
v_1)}$, as defined on page 29 of \cite{KL}. It is holomorphic, and 
we have

{\bf Theorem:} If $n \in \Z_{> 0}$ is such that
\[
v \in \left( \frac{1}{n} \Z \right)^2 - \Z^2 \; ,
\]
then the $(12n)$-th power of $\Si^v$ is a non-vanishing algebraic 
function on $Y (n)_{\C}$. It descends to $Y (n)$, viewed as a geometrically 
connected scheme over the subfield $\Q \left( e^{\frac{2\pi i}{n}} 
\right)$ of $\C$:
\[
\Si^{v} \in \Gamma (Y (n) , \Oh^{\ast}_{Y (n)}) \otimes_{\Z} \Z 
\left[ \frac{1}{12n} \right] \; .
\]
As an element of $\Gamma (Y (n) , \Oh^{\ast}_{Y (n)}) \otimes_{\Z} 
\Z \left[ \frac{1}{12n} \right]$, the function $\Si^{v}$ only 
depends on $v \mod \Z^2$.

{\bf Proof:} The first statement is \cite{KL}, II, Theorem 1.2 -- 
but note from Theorem {\rm \ref{41}}, we know that {\it some} power 
of $\Si^{v}$ is an algebraic function on $Y (n)_{\C}$. For the 
descent to $Y (n)$, one uses the $q$-expansion principle. The last 
statement follows from Lemma {\rm \ref{13}} -- again, the 
independence of $\Si^{v}$ in
\[
\Gamma (Y (n) , \Oh^{\ast}_{Y (n)}) \otimes_{\Z} \Q
\]
is also predicted by Theorem {\rm \ref{41}}. \beweisende
\end{punkt}

%
%
\begin{punkt}\label{45} \em Let $v \in \left( \frac{1}{n} \Z \right)^2 - 
\Z^2$, and consider the section
\begin{eqnarray*}
\Hh^+ & \longrightarrow & \C \times \Hh^+ \; ,\\
\tau & \longmapsto & (-v_2 \tau + v_1 , \tau) \; .
\end{eqnarray*}
It descends to the level of $Y (n)_{\C}$ and $\Eh_{n,\C}$, defining 
a non-zero $n$-torsion section $i_{v,\C}$ of $\tau_{n,\C}$. Via the 
canonical embedding of $\Q \left( e^{\frac{2\pi i}{n}} \right)$ 
into $\C$, we get
\[
i_v : Y (n) \longrightarrow \Eh_n \; .
\]
By {\rm \ref{41}} and \cite{KL}, II, Proposition 1.3, the element
\[
\varphi (\{ i_v \} ) \in \Ext^1_{HDR^s_{\Q} (Y (n)_{\C})} (\Q (0) , 
\Q (1))
\]
satisfies: the underlying collection
\[
\For (\varphi (\{ i_v \})) \quad \in \prod_{\sigma : \Q ( e^{\frac{2\pi 
i}{n}}) \hookrightarrow \C} \Ext^1_{\Var_{\Q} (Y (n)_{\C})} 
(\Q (0) , \Q (1))
\]
of extensions of variations underlying $\varphi (\{ i_v \} )$ lies 
in the image of
\[
\Delta : \Gamma (Y (n) , \Oh^{\ast}_{Y (n)}) \otimes_{\Z} \Q 
\longrightarrow \prod_{\sigma} \Gamma (Y (n) , \Oh^{\ast}_{Y (n)} ) 
\otimes_{\Z} \Q \; .
\]
From Proposition {\rm \ref{23}}, we conclude:

{\bf Proposition:} For any $v \in \left( \frac{1}{n} \Z \right)^2 - 
\Z^2$, we have
\[
\varphi (\{ i_v \}) \in \Gamma (Y (n) , \Oh^{\ast}_{Y (n)}) 
\otimes_{\Z} \Q \; .
\]
Via the canonical embedding of $\Q \left( e^{\frac{2\pi i}{n}} 
\right)$ into $\C$, we have the equality
\[
\varphi ( \{ i_v \}) = (\Si^v)^{-1} \; .
\]
\end{punkt}

\begin{punkt}\label{46} \em We now work over an arbitrary base 
field $k$ which is embeddable into $\C$, but still assume that for 
some $n \ge 3$, the whole $n$-torsion of
\[
\pi : \Eh \longrightarrow B
\]
consists of sections. Let $s \in \widetilde{\Eh} (B)$ be an $n$-torsion section. It comes 
about as the base change by some $f: B \longrightarrow Y(n)$ of a section
\[
i_v : Y (n) \longrightarrow \Eh_n
\]
of $\pi_n$. Because of the functoriality statement in Theorem {\rm 
\ref{33}}, we know that
\[
\varphi (\{ s \}) = f^{\ast} \varphi (\{ i_v \}) \; .
\]
$\varphi (\{ s \})$ is therefore an algebraic function on $B$.
\end{punkt}

\begin{punkt}\label{47} \em We now remove the hypothesis on the
$n$-torsion in {\rm \ref{46}}. So let $\Eh$ be an arbitrary 
elliptic curve over $B$, and $s \in \widetilde{\Eh} (B)$ a torsion 
section. Choose a multiple $n \ge 3$ of the order of $s$, and a 
finite \'etale Galois 
covering $C$ of $B$ such that the whole $n$-torsion of $\Eh 
\times_B C$ consists of sections. For the base change $s_C$ of $s$ 
to $C$, we have, by {\rm \ref{46}}
\[
\varphi (\{ s_C \}) \in \Gamma (C , \Oh^{\ast}_C) \otimes_{\Z} \Q \; .
\]
Using the functoriality of $\varphi$ with respect to base change 
under the automorphisms of $C$ over $B$, one concludes purely 
formally:

{\bf Theorem:} Let $s \in \widetilde{\Eh} (B)$ a torsion section. Then 
\[
\varphi (\{ s \}) \in \Ext^1_{HDR^s_{\Q} (B)} (\Q (0) , \Q (1))
\]
lies in the image of $\kappa_B$:
\[
\varphi (\{ s \}) \in \Gamma (B , \Oh^{\ast}_B) \otimes_{\Z} \Q \; .
\]
\end{punkt}

\begin{punkt}\label{48} \em For a scheme $B$, which is smooth, separated,
connected and of finite type over some field of characteristic $0$, 
fix a geometric point $\overline{b}$, and consider the projective 
system $\{ B_{\alpha} \tei \alpha \in I \}$ of pointed finite \'etale 
coverings of $B$. We have
\[
\Gamma (B , \Oh^{\ast}_B) \otimes_{\Z} \Q = \left( \limr{\alpha} 
\Gamma (B_{\alpha} , \Oh^{\ast}_{B_{\alpha}}) \otimes_{\Z} \Q 
\right)^{\pi_1 (B , \overline{b})} \; .
\]
Let $\Eh$ be an elliptic curve over $B$. We denote by $\Lh_0 (\Eh)$ 
the subspace of $\ker (d)$ of divisors supported on torsion sections.

{\bf Definition:} The subspace of {\it elliptic modular units of 
the zeroeth kind on $B$} is defined as
\[
\left( \limr{\alpha} \varphi (\Lh_0 (\Eh \times_B B_{\alpha})) 
\right)^{\pi_1 (B , \overline{b})} \subset \Gamma (B , 
\Oh^{\ast}_B) \otimes_{\Z} \Q \; .
\]
It is denoted by $EM_0 (\Eh)$.
\end{punkt}

\begin{punkt}\label{49} \em It is natural to ask for the size of 
$EM_0 (\Eh)$ inside $\Gamma (B , \Oh^{\ast}_B) \otimes_{\Z} \Q$.

{\bf Examples:} a) For $\Eh_n , n \ge 3$, the group
\[
EM_0 (\Eh_n) \subset \Gamma (Y (n), \Oh^{\ast}_{Y (N)}) 
\otimes_{\Z} \Q
\]
contains the Siegel units $\Si^v$. By \cite{KL}, IV, Theorem 1.1, 
we thus have
\[
\Q \left( e^{\frac{2\pi i}{n}} \right) \cdot EM_0 (\Eh_n) = \Gamma 
(Y (n) , \Oh^{\ast}_{Y (n)}) \otimes_{\Z} \Q \; .
\]
The norm compatibility statement {\rm \ref{32}} b) translates into 
a distribution relation modulo roots of unity of the Siegel units. 
For the precise distribution relation, see \cite{Ku}, Theorem 2.2.

b) Let $\Eh = E$ be an elliptic curve over a number field $F$ with 
complex multiplication by $K \subset F$, such that every torsion 
point of $E$ is defined over the maximal abelian extension 
$K^{\ab}$ of $K$.
Then $EM_0 (E \otimes_F F^{\ab})$ contains the classical {\it 
elliptic units} modulo torsion (see \cite{Sh}, II, \S\,2).
\end{punkt}

\newpage
\section{Elliptic modular units of the first kind}

\begin{punkt}\label{51} \em Suppose given a section $s \in 
\widetilde{\Eh} (B)$, and a divisor
\[
D = \sum_t q_t \cdot (t) \in \Div^0 (\Eh) \otimes_{\Z} \Q \; ,
\]
where $t$ runs through the torsion sections of
$\pi : \Eh \longrightarrow B$.
Assume that $q_t = 0$ if $s -t$ is not disjoint from the zero 
section. Then
\[
S_{D,s} := \sum_t q_t \{ s-t \} \in \Lh (\Eh)
\]
is an element of the kernel of $d$.

{\bf Theorem:} $\varphi (S_{D,s}) \in \Ext^1_{HDR^s_{\Q} (B)} (\Q 
(0) , \Q (1))$ lies in the image of $\kappa_B$:
\[
\varphi (S_{D,s}) \in \Gamma (B , \Oh^{\ast}_B) \otimes_{\Z} \Q \; .
\]

{\bf Proof:} Our situation arises via the base change
$s : B \longrightarrow \widetilde{\Eh}$
from the projection
\[
\pr_1 : \widetilde{\Eh} \times_B \Eh \longrightarrow \widetilde{\Eh} \; ,
\]
with the torsion sections $t$ replaced by
\begin{eqnarray*}
t : \widetilde{\Eh} & \longrightarrow & \widetilde{\Eh} \times_B \Eh \; , \\
x & \longmapsto & (x,t) \; ,
\end{eqnarray*}
and $s$ replaced by $\Delta$. So we need to show the statement for 
$\pr_1$ and $S_{\pr^{\ast}_1 D, \Delta}$.

As in the proof of \ref{47}, we may 
assume that there is a torsion section $t'$
not contained in the support of $D$. Then base change of
\[
\pr_1 : \widetilde{\Eh} \times_B \Eh \longrightarrow \widetilde{\Eh}
\]
via $t'$ gives back the original situation, with $s$ replaced by 
$t'$.

Our claim then follows from {\rm \ref{47}} and {\rm \ref{24}}. \beweisende
\end{punkt}

\begin{punkt}\label{52} \em We want to write down explicit formulae 
for the $\varphi (S_{D,s})$. In section 4, our geometrical object 
of study was the moduli space $Y (n)$ for $n$-torsion sections. The 
moduli space for $n$-torsion sections plus an additional section is 
the universal elliptic curve $\Eh_n$ itself.

So let $n \ge 3$, and
\[
\pr_1 : \Eh_n \times_{Y (n)} \Eh_n \longrightarrow \Eh_n \; .
\]
Recall from {\rm \ref{45}} that we parameterized the $n$-torsion 
sections of $\pi_n$ by $v \in \left( \frac{1}{n} \Z \right)^2 / \Z^2$:
\[
i_v : Y (n) \longrightarrow \Eh_n \; .
\]
Via base change, we get $n$-torsion sections $i_v$ of $\pr_1$.

Assume given a divisor
\[
D = \sum_v q_v \cdot (i_v)
\]
of degree $0$. Then
\[
S_D := S_{D,\Delta} = \sum_v q_v \cdot (\Delta - i_v) \in \Lh (U_D 
\times_{Y (n)} \Eh_n) \; ,
\]
where $U_D \subset \Eh_n$ is the open subscheme of $\Eh_n$ 
complementary to the support of $D$. We have $S_D \in \ker (d)$.

We apply base change to $\C$ and determine
\[
\varphi (S_D)^{MHS} \in \Gamma (U_{D,\C} , \Oh^{\ast}_{D,\C}) 
\otimes_{\Z} \Q \; .
\]
By {\rm \ref{312}}, it differs from the inverse of the multivalued 
function
\[
\Si_D : (z,\tau) \longmapsto \prod_v \Si (z + v_2 \tau - v_1 , 
\tau)^{q_v}
\]
by a multivalued function of the shape
$\exp (2\pi i \cdot F)$,
where $F$ is a polynomial in $\Q [r_1 , r_2]$ of total degree 
smaller or equal to one.

{\bf Definition:} A multivalued function is called a {\it 
holomorphic modification} of $\Si_D$, if it is of the shape
\[
\Si_D \cdot \exp (2\pi i \cdot F)
\]
for a polynomial $F$ in $\Q [r_1 , r_2]$ of total degree smaller or 
equal to one.

Since expressions of the form $\exp (2\pi i \cdot F)$ are 
holomorphic on $\C \times \Hh^+$ if and only if $F = 0$, we have:

{\bf Proposition:} There exists a unique holomorphic modification 
of $\Si_D$. It is equal to the inverse of $\varphi (S_D)^{MHS}$.
\end{punkt}

In practical terms, this means: in order to find $(\varphi 
(S_D)^{MHS})^{-1}$, write down the formula for $S_D$ and cancel 
the factors
\[
\exp (2\pi i \cdot (u_1 r_1 + u_2 r_2))
\]
of non-holomorphicity.

\begin{punkt}\label{57} \em Now let again $\Eh$ be an elliptic 
curve over a $k$-scheme $B$, $k$ being of characteristic $0$. We 
define $\Lh_1 (\Eh)$ as the subspace of $\ker (d)$ generated by 
divisors of the shape $S_{D,s}$ as in {\rm \ref{51}}.

{\bf Definition:} The subspace of {\it elliptic modular units of 
the first kind on $B$} is defined as
\[
\left( \limr{\alpha} \varphi (\Lh_1 (\Eh \times_{B} B_{\alpha} )) 
\right)^{\pi_1 (B ,\overline{b})} \subset \Gamma (B , \Oh^{\ast}_B) 
\otimes_{\Z} \Q \; .
\]
It is denoted by $EM_1 (\Eh)$.

{\bf Example:} Let $s \in \widetilde{\Eh} (B)$ and $N \ge 1$ such that 
$[N] s$ is still in $\widetilde{\Eh} (B)$, i.e., disjoint from the zero 
section. Then
\[
S := \{ [N] s \} - N^2 \{ s \}
\]
is in $\ker (d)$. We claim that $\varphi (S)$ is actually an elliptic 
modular unit of the first kind.

In order to see this, we may assume that the whole $N$-torsion of 
$\Eh$ consists of sections of $\pi$. Then
$T := \{ [N] s\} - \sum_{t \in \Eh [N] (B)} \{ s- t \}$
is in $\ker (d)$. So our claim follows if we show {\it strong norm 
compatibility:}
\[
\varphi \left( T 
\right) = 1 \in \Gamma (B , \Oh^{\ast}_B) \otimes_{\Z} \Q \; .
\]
This is achieved, as in the proof of Theorem {\rm \ref{51}}, by
reducing to the case when $s$ is torsion.
\end{punkt}
\newpage
\section{Elliptic modular units of the second kind}

\begin{punkt}\label{61} \em The last special elements of $\ker (d)$ 
we want to consider are the parallelograms
\[
\{ s+t \} + \{ s-t \} - 2\{ s \} - 2 \{ t \} \; ,
\]
for $s,t, s+t , s-t \in \widetilde{\Eh} (B)$.

{\bf Theorem:} Let $S_P := \{ s+t \} + \{ s-t \} - 2\{ s \} - 2\{ t 
\}$, for $s$ and $t$ as above. Then
\[
\varphi (S_P) \in \Ext^1_{HDR^s_{\Q} (B)} (\Q (0) , \Q (1))
\]
lies in the image of $\kappa_B$:
\[
\varphi (S_P) \in \Gamma (B , \Oh^{\ast}_B) \otimes_{\Z} \Q \; .
\]

{\bf Proof:} If $s$ and $t$ are torsion sections, then the claim 
holds by Theorem {\rm \ref{47}}. The general case follows as in 
Theorem {\rm \ref{51}}. We leave the details to the reader. 
\beweisende
\end{punkt}

\begin{punkt}\label{62} \em Again, we write down an explicit 
formula for $\varphi (S)$. Following the procedure of {\rm 
\ref{52}}, we work on the moduli space for $n$-torsion sections, 
plus two additional sections.

Let $n \ge 3$, and
\[
\pr_{12} : \Eh_n \times_{Y (n)} \Eh_n \times_{Y (n)} \Eh_n 
\longrightarrow \Eh_n \times_{Y (n)} \Eh_n \; .
\]
There are two sections $s$ and $t$ of $\pr_{12}$:
\begin{eqnarray*}
s : (x,y) & \longmapsto & (x,y,x) \; , \\
t : (x,y) & \longmapsto & (x,y,y) \; .
\end{eqnarray*}
Write
\[
S_P := \{ s+t \} + \{ s-t \} - 2\{ s \} - 2\{ t \} \; .
\]
It is an element of $\Lh (V \times_{Y (n)} , \Eh_n)$, where $V 
\subset \Eh_n \times_{Y (n)} \Eh_n$ is the complement of the union 
of the two zero sections, of the diagonal $\Delta$, and of the 
anti-diagonal
\[
\Delta^a := \{ (x , -x) \in \Eh_n \times_{Y (n)} \Eh_n \tei x \in 
\Eh_n \} \; .
\]
We apply base change to $\C$ and determine
\[
\varphi (S_P)^{MHS} \in \Gamma (V , \Oh^{\ast}_V) \otimes_{\Z} \Q \; .
\]

In the following, write
\[
f : X \sto Y
\]
for a morphism of real manifolds if it is clear from the context 
which closed subset of $X$ one would have to remove in order to get 
the domain of definition for $f$.

{\bf Definition:} Define
\begin{eqnarray*}
\Si_P : \C \times \C \times \Hh^+ & \sto & \C \; , \\
(z_1 , z_2 , \tau) & \longmapsto & \frac{\Si (z_1 + z_2 , \tau) 
(\Si (z_1 - z_2 , \tau)}{\Si (z_1 , \tau)^2 \Si (z_2,\tau)^2 } \; .
\end{eqnarray*}

{\bf Proposition:} $\Si_P$ is holomorphic on its domain of 
definition. It is equal to the inverse of $\varphi (S_P)^{MHS}$.

{\bf Proof:} The proof of the first claim uses Lemma {\rm 
\ref{13}}. The second claim follows from {\rm \ref{312}}. \beweisende
\end{punkt}

\begin{punkt}\label{64} \em Let $B$ be connected, smooth and separated 
of finite type over some field of characteristic $0$,
\[
\pi : \Eh \longrightarrow B
\]
an elliptic curve. We define $\Lh_2 (\Eh)$ as the subspace of $\ker 
(d)$ generated by divisors of the shape
\[
\{ s + t \} + \{ s-t \} - 2\{s \} - 2 \{ t \} \; ,
\]
for $s,t , s+t , s-t \in \widetilde{\Eh} (B)$.

Again, write $\{ B_{\alpha} \tei \alpha \in I \}$ for the 
projective system of pointed finite \'etale coverings of $B$.

{\bf Definition:} The subspace of {\it elliptic modular units of 
the second kind on $B$} is defined as
\[
\left( \limr{\alpha} \varphi (\Lh_2 (\Eh \times_B B_{\alpha})) 
\right)^{\pi_1 (B , \overline{b})} \subset \Gamma (B , 
\Oh^{\ast}_B) \otimes_{\Z} \Q \; .
\]
It is denoted by $EM_2 (\Eh)$.
\end{punkt}
\newpage
\section{Proof of Theorem \ref{32}}

\begin{punkt}\label{72} \em Let $A$ be an abelian group
\[
\delta_A : \Q [A]  \longrightarrow  A \otimes_{\Z} A \otimes_{\Z} 
\Q \; , \quad
\{ a \}  \longmapsto  a \otimes a \; ,
\]
and denote by $d_A$ the restriction of $\delta_A$ to 
$\Lh (A) := \Q [A - \{ 0 \} ]$.

{\bf Proposition:} 
The kernel of $d_A$ is generated by expressions of the 
following shape:

o) $\{ a \}$ for $a \in A_{\tors}$, $a \neq 0$.

i) $\{ a \} - \{ a-b \}$ for $a \in A$, $b \in A_{\tors}$, $a, a-b \neq 0$;

\quad $\{ Na \} - N^2 \{ a \}$ for $a \in A$, $N \ge 1$, $Na \neq 0$.

ii) $\{ a+b \} + \{ a-b \} - 2\{ a \} - 2\{ b \}$ for $a,b \in A$, 
$a,b, a+b , a-b \neq 0$.

{\bf Proof:} First show that the set of expressions as in o)--ii),
with the condition ``$\ne 0$'' removed, generates $\ker (\delta_A)$.
Then observe that
$\delta_A$
factors through the projection
\[
\Q[A] \longrightarrow \Q[A - \{ 0 \} ]
\]
given by sending $\{ 0 \}$ to $0$.                 \beweisende
\end{punkt}

\begin{punkt}\label{73} \em We are finally able to show:

{\bf Theorem:} The morphism
\[
\varphi = \varphi (\Eh) : \ker (d) \longrightarrow 
\Ext^1_{HDR^s_{\Q} (B)} (\Q (0) , \Q (1))
\]
factors through $\Gamma (B , \Oh^{\ast}_B) \otimes_{\Z} \Q$.

{\bf Proof:} We may assume that $B$ is connected. This implies that 
if $S$ is a non-zero section of
\[
\pi : \Eh \longrightarrow B \; ,\]
then it is actually disjoint from the zero section $i$ on an open 
dense subscheme of $B$.

So if $S \in \ker (d)$, Proposition {\rm \ref{72}}, and {\rm \ref{47}}, {\rm \ref{51}}, {\rm \ref{57}},
and {\rm \ref{61}} tell us that the restriction of $\varphi (S)$ to some 
open dense subscheme $B'$ of $B$ lies in
\[
\Gamma (B' , \Oh^{\ast}_{B'}) \otimes_{\Z} \Q \; .
\]
Our claim follows from Theorem {\rm \ref{24}}. \beweisende
\end{punkt}
We thus get the proof of parts a) and b) of Theorem {\rm \ref{32}}.

\begin{punkt}\label{74} \em Using the explicit formulae of {\rm \ref{42}}, 
{\rm \ref{52}}, and {\rm \ref{62}}, and the classical relations between
the Siegel function and the Weierstra{\ss} function and its 
derivatives, one identifies $\varphi$ on the special elements
$\{s , t\} \in \ker (d)$ of \cite{GL}, 4.5 (cf.\ e.g.\ \cite{R}, 
V.4, proof of Proposition 4.4). One gets a relative version
of \cite{W3}, Proposition 1.9.1, which actually shows that the present
construction and that of loc.\ cit.\ produce the same functions. 

In order to show 
the relation to the non-archimedian local heights in {\rm \ref{32}}.c),
one then proceeds as in loc.\ cit.
\end{punkt}

\begin{punkt}\label{77} \em Let $B$ be a connected scheme, which is
smooth and separated 
of finite type over some field of characteristic $0$,
\[
\pi : \Eh \longrightarrow B
\]
an elliptic curve. Write $\{ B_{\alpha} \tei \alpha \in I \}$ for 
the projective system of pointed finite \'etale coverings of $B$.

{\bf Definition:} The subspace of {\it elliptic modular units on 
$B$} is defined as
\[
\left( \limr{a} \varphi (\ker (d (\Eh \times_B B_{\alpha}))) 
\right)^{\pi_1 (B, \overline{b})} \subset \Gamma (B , \Oh^{\ast}_B) 
\otimes_{\Z} \Q \; .
\]
It is denoted by $EM (\Eh)$.
\end{punkt}

\begin{punkt}\label{78} \em In the case when $B$ is a point, the
following observation is due to Goncharov and Levin (\cite{GL}, Corollary
4.5):

{\bf Theorem:}
Zariski-locally on $B$, every non-vanishing function  
is an elliptic modular unit. 

{\bf Proof:} Assume given $s$, $t \in \widetilde{\Eh}$ such that $s$ is
disjoint from both $t$ and $-t$. Then we have
\[
\{s , t\} = \{ s+t \} + \{ s-t \} - 2\{ s \} - 2\{ t \} \in \ker(d) \; .
\]
According to the relative version of \cite{W3}, Proposition 1.9.1,
\[
\varphi (\{ s,t \}) = ((x(t) - x(s)) \Delta^{-1/6})^{-1}
\]
for any local Weierstra{\ss} equation
$y^2 = x^3 + ax + b$
of $\Eh$.                             \beweisende
\end{punkt}

\end{document}